\documentclass{amsart}
\usepackage{cases}
\usepackage{amsmath}
\usepackage{amssymb}
\usepackage{amsfonts}
\usepackage{graphicx}
\newtheorem{theorem}{Theorem}[section]

\newtheorem{corollary}[theorem]{Corollary}

\newtheorem{definition}[theorem]{Definition}
\newtheorem{example}[theorem]{Example}

\newtheorem{lemma}[theorem]{Lemma}

\newtheorem*{problem*}{Problem 1}
\newtheorem{proposition}[theorem]{Proposition}
\newtheorem{remark}[theorem]{Remark}

\newtheorem{thm}{Theorem}

\begin{document}
\title[Splitting theorems with
Bakry-\'{E}mery curvature]
{Splitting theorems on complete manifolds with
Bakry-\'{E}mery curvature}
\author{Jia-Yong Wu}

\address{Department of Mathematics, Shanghai Maritime University,
Haigang Avenue 1550, Shanghai 201306, P. R. China}

\email{jywu81@yahoo.com}

\thanks{This work is partially supported by the NSFC11101267.}


\dedicatory{}

\begin{abstract}
In this paper we study some splitting properties on complete noncompact
manifolds with smooth measures when $\infty$-dimensional
Bakry-\'{E}mery Ricci curvature is bounded from below by some negative
constant and spectrum of the weighted Laplacian has a positive lower
bound. These results extend the cases of Ricci curvature and
$m$-dimensional Bakry-\'{E}mery Ricci curvature.
\end{abstract}
\maketitle

\section{Introduction and main results}
Let $(M,g)$ be an $n$-dimensional complete Riemannian manifold
and $\varphi$ be a smooth function on $M$. We define
the weighted Laplacian
\[
\Delta_\varphi:=\Delta-\nabla \varphi\cdot\nabla,
\]
which is the infinitesimal generator of the Dirichlet form
\[
\mathcal {E}(\phi_1,\phi_2)=\int_M \langle\nabla \phi_1, \nabla \phi_2\rangle d\mu,
\,\,\,\forall \phi_1, \phi_2\in C_0^{\infty}(M),
\]
where $\mu$ is an invariant measure of $\Delta_\varphi$ given by $
d\mu=e^{-\varphi}dv(g)$. Clearly, the weighted Laplacian is
self-adjoint with respect to the weighted measure $d\mu$ and
$(M,g,e^{-\varphi}dv)$ is a smooth metric measure spaces.
A smooth function $u$ is called weighted harmonic if
$\Delta_\varphi u=0$. The first nontrivial eigenvalue of the
weighted Laplacian in $(M,g,e^{-\varphi}dv)$ is defined by
\[
\lambda_1(M):={\inf\limits_{f\neq 0}}\left\{{\mathcal {E}(f,f):
\int_Mf^2d\mu=1,\int_Mfd\mu=0}
\right\}.
\]
The above infimum can be achieved by some smooth
eigenfunction $f$, which satisfies the following Euler-Lagrange
equation
\[
\Delta_\varphi f=-\lambda_1f.
\]

The $m$-dimensional Bakry-\'{E}mery Ricci curvature
(see also \cite{[BE],[BQ1],[BQ2],[LD]}) is linked with
the smooth metric measure space $(M,g,e^{-\varphi}dv)$, which
is defined by
\[
Ric_{m,n}:=Ric +Hess(\varphi)-\frac{\nabla \varphi \otimes \nabla
\varphi}{m-n},
\]
where $Ric$ and $Hess$ denote the Ricci curvature and the Hessian of
the metric $g$, respectively. Here $m\geq n$ is a constant, and
$m=n$ if and only if $\varphi$ is a constant \cite{[LD],[Lott]}.
If we let $m$ be infinite, then the $m$-dimensional Bakry-\'{E}mery Ricci
curvature becomes the $\infty$-dimensional Bakry-\'{E}mery Ricci
curvature
\[
Ric_\infty:=\lim_{m\to \infty}Ric_{m,n}=Ric+Hess(\varphi).
\]
This curvature is closely related to the gradient Ricci soliton
\[
Ric_\infty=\lambda g
\]
for some constant $\lambda$. The soliton is called expanding, steady
and shrinking, accordingly, if $\lambda<0$, $\lambda=0$ and
$\lambda>0$, which plays an important role in the theory of
Ricci flows \cite{[Cao1]}. The gradient Ricci solitons are
special solutions of Ricci flows and often arise from the
blow up analysis of the singularities of Ricci flows
\cite{[Hamilton]}.

If $\varphi$ is constant, the above weighted geometric quantities
all return to the classical case. Taking this point of view,
it is natural to ask if classical results involving Ricci
curvature in the geometric analysis remain valid in the
Bakry-\'{E}mery curvature case. Along this direction, a lot
of excellent work actually obtained in the past several years.
Interesting generalized results concerning Cheeger-Gromoll splitting
theorems \cite{[Che-Gro1],[Che-Gro2]} are Cheeger-Gromoll type
splitting theorems with $Ric_\infty\geq 0$ showed by Wei and
Wylie \cite{WeiWy} that is originally due to Lichnerowicz.
This result was then improved by Fang, Li and Zhang \cite{[FLZ]}.
Recently Munteanu and Wang \cite{[MuWa],[MuWa2]}
studied function theoretic and spectral properties on complete
noncompact smooth metric measure space with $\infty$-dimensional
Bakry-\'{E}mery Ricci curvature bounded below. In particular,
they obtained some splitting results on complete
noncompact gradient steady and expanding Ricci solitons.

In \cite{[Wang2]}, Wang proved a splitting theorem on complete
metric measure spaces with $m$-dimensional Bakry-\'{E}mery
Ricci curvature bounded below by a negative multiple of
the lower bound of the weighted spectrum. In particular,
he proved that

\begin{thm}\label{Wmain}
Let $(M^n,g,e^{-\varphi}dv)$ be a complete smooth metric measure space
of dimension $n\geq3$ with $Ric_{m,n}\geq-(m-1)$. If
$\lambda_1(M)\geq(m-2)$, then either
\begin{enumerate}
  \item $M$ has only one end with infinite weighted volume; or

\item $M=\mathbb{R}\times N$ with the warped product metric
\[
ds_M^2=dt^2+\cosh^2tds_N^2,
\]
where $N$ is an $(n-1)$-dimensional compact Riemannian manifold.
In this case, $\lambda_1(M)=m-2$.
\end{enumerate}
\end{thm}
This result extended Li-Wang's theorem \cite{[Li-Wang1]}
to the weighted measure case. If $\lambda_1(M)>(m-2)$, Theorem
\ref{Wmain} asserts that $M$ has only one end with infinite
weighted volume. Naturally, we may ask if finite weighted
volume ends can be ruled out in this case. In \cite{[Wu2]},
the author proved that there exists an example to indicate that 
finite weighted volume ends can exist in Theorem \ref{Wmain}. 
Moreover this example seems to be the only case when $M$
has a finite weighted volume end if the weighted eigenvalue
has an optimal positive lower bound and $m$-dimensional
Bakry-\'{E}mery Ricci curvature is bounded below by
some negative constant. Precisely,
\begin{thm}\label{wth}
Let $(M^n,g,e^{-\varphi}dv)$ be a complete smooth metric measure space
of dimension $n\geq3$ with $Ric_{m,n}\geq-(m-1)$. If
$\lambda_1(M)\geq\frac{(m-1)^2}{4}$, then either
\begin{enumerate}
  \item $M$ has only one end; or

\item $M=\mathbb{R}\times N$ with the warped product metric
\[
ds_M^2=dt^2+\exp(-2t)ds_N^2
\]
where $N$ is an $(n-1)$-dimensional compact manifold. Moreover,
\[
\varphi(t,x)=\varphi(0,x)+(m-n)t
\]
for all $(t,x)\in\mathbb{R}\times N$.
\end{enumerate}
\end{thm}
Theorem \ref{wth} was also independently proved by Su and Zhang
\cite{[SuZh]}. This result can also be viewed as a weighted version
of Li-Wang's theorem \cite{[Li-Wang2]}. In \cite{[Wu]}
(see also \cite{[Wang1]}) we know that if $Ric_{m,n}\geq-(m-1)$,
then $\lambda_1(M)\leq\frac{(m-1)^2}{4}$.
Hence the eigenvalue assumption in Theorem \ref{wth} is,
in fact that $\lambda_1(M)=\frac{(m-1)^2}{4}$.

\vspace{0.5em}

In this paper, we continue our study of splitting type theorems
for a smooth metric measure space, now under slight different assumption
that its $\infty$-dimensional Bakry-\'{E}mery Ricci curvature is bounded
from below. The main result is
\begin{theorem}\label{Rmain}
Let $(M^n,g,e^{-\varphi}dv)$ be a complete smooth metric measure space
of dimension $n\geq3$ with $Ric_\infty\geq-(n-1)$. Assume that
$|\nabla\varphi|\leq\theta$ for some positive constant $\theta$.
For any $k>n$, if
\[
\lambda_1(M)\geq\frac{k-2}{k-1}
\left(\frac{\theta^2}{k-n}+n-1\right),
\]
then either
\begin{enumerate}
  \item $M$ has only one end with infinite weighted volume; or

\item $M=\mathbb{R}\times N$ with the warped product metric
\begin{equation}\label{biaoxing}
ds_M^2=dt^2+\cosh^2\left[\sqrt{\frac{\theta^2}{(k-1)(k-n)}
+\frac{n-1}{k-1}}\,\,t\right]ds_N^2,
\end{equation}
where $N$ is an $(n-1)$-dimensional compact Riemannian manifold.
In this case, $\lambda_1(M)=\frac{k-2}{k-1}
\left(\frac{\theta^2}{k-n}+n-1\right)$.
\end{enumerate}
\end{theorem}
In Theorem \ref{Rmain}, if $\theta=0$, by letting $k\to n$, we then have
Li-Wang's theorem \cite{[Li-Wang1]}. In particular, if we let
$k=n+\theta$, then we have a explicit version.
\begin{corollary}\label{Rcor}
Let $(M^n,g,e^{-\varphi}dv)$ be a complete smooth metric measure space
of dimension $n\geq3$ with $Ric_\infty\geq-(n-1)$. Assume that
$|\nabla\varphi|\leq\theta$ for some positive constant $\theta$.
If $\lambda_1(M)\geq n-2+\theta$, then either
\begin{enumerate}
  \item $M$ has only one end with infinite weighted volume; or

\item $M=\mathbb{R}\times N$ with the warped product metric
\begin{equation}\label{biaon2}
ds_M^2=dt^2+\cosh^2tds_N^2,
\end{equation}
where $N$ is an $(n-1)$-dimensional compact Riemannian manifold.
In this case, $\lambda_1(M)=n-2+\theta$.
\end{enumerate}
\end{corollary}
A natural question as above discussions is that if finite
weighted volume ends can be excluded in Corollary \ref{Rcor}.
Similar to Theorem \ref{wth}, we assert that finite weighted
volume ends in Corollary \ref{Rcor} may occur by the following example.
\begin{example}
Consider the $n$-dimensional complete manifold
$M=\mathbb{R}\times N$ endowed with the warped product metric
\[
ds_M^2=dt^2+\exp(-2t)ds_N^2.
\]
If $\{\bar{e}_\alpha\}$ for $\alpha=2,...,n$ form an orthonormal basis of
the tangent space of $N$, then $e_1=\frac{\partial}{\partial t}$ together
with $\{e_\alpha=\exp(-t)\bar{e}_\alpha\}$ form an orthonormal basis for
the tangent space of $M$. By the standard computation,
\[
Ric_{M,1j}=-(n-1)\delta_{1j}
\]
and
\[
Ric_{M,\alpha\beta}=\exp(2t)Ric_{N,\alpha\beta}-(n-1)\delta_{\alpha\beta}.
\]
If we let the weighted function $\varphi:=\theta t$, then the
$\infty$-dimensional Bakry-\'{E}mery Ricci curvature of $M$ is
\[
Ric_{\infty,1j}=Ric_{M,1j}+\varphi_{1j}
=-(n-1)\delta_{1j}
\]
and
\[
Ric_{\infty,\alpha\beta}=\exp(2t)Ric_{N,\alpha\beta}-(n-1)\delta_{\alpha\beta}.
\]
Clearly, $Ric_N\geq 0$ if and only if $Ric_{\infty}\geq-(n-1)$.
Moreover, we \emph{claim} that
\[
\lambda_1(M)=\frac{(n-1+\theta)^2}{4}.
\]
Indeed, we may choose the function
\[
f=\exp\left(\frac{n-1+\theta}{2}t\right).
\]
A direct computation yields that
\begin{equation*}
\begin{aligned}
\Delta_\varphi f&=\frac{d^2f}{dt^2}-(n-1)\frac{df}{dt}
-\frac{d\varphi}{dt}\cdot\frac{df}{dt}\\
&=-\frac{(n-1+\theta)^2}{4}f,
\end{aligned}
\end{equation*}
since $\Delta=\frac{\partial^2}{\partial t^2}
-(n-1)\frac{\partial}{\partial t}+\exp(2t)\Delta_N$.
So we have $\lambda_1(M)\geq\frac{(n-1+\theta)^2}{4}$
by Proposition 1.4 of \cite{[Wu1]}. Combining this and
Theorem \ref{wthm}, we immediately
conclude that $\lambda_1(M)=\frac{(n-1+\theta)^2}{4}$ as claimed.
\end{example}

We can show that the above result is the only case (possible with
different weighted function) when $M$ has a finite weighted volume
end if $M$ achieves an optimal positive lower bound of the
weighted spectrum.
\begin{theorem}\label{wth2}
Let $(M^n,g,e^{-\varphi}dv)$ be a complete smooth metric measure space
of dimension $n\geq3$ with $Ric_\infty\geq-(n-1)$. Assume that
$|\nabla\varphi|\leq\theta$ for some positive constant $\theta$.
If $\lambda_1(M)\geq\frac{(n-1+\theta)^2}{4}$, then either
\begin{enumerate}
  \item $M$ has only one end; or

\item $M=\mathbb{R}\times N$ with the warped product metric
\[
ds_M^2=dt^2+\exp(-2t)ds_N^2,
\]
where $N$ is an $(n-1)$-dimensional compact manifold with
nonnegative Ricci curvature. Moreover,
\[
\varphi(t,x)=\varphi(0,x)+\theta t
\]
for all $(t,x)\in\mathbb{R}\times N$.
\end{enumerate}
\end{theorem}

We remark that the assumption on $\lambda_1(M)$ in Theorem \ref{wth2}
in fact only occurs on the equality case due to the following result
(see Theorem 1.1 in \cite{[Wu1]}), which is a generalization of Cheng's
eigenvalue estimate \cite{[Cheng]}.
\begin{thm}\label{wthm}
Let $(M^n,g,e^{-\varphi}dv)$ be a complete smooth metric measure space
of dimension $n\geq2$ with $Ric_\infty\geq-(n-1)$. Assume that
$|\nabla\varphi|\leq\theta$ for some positive constant $\theta$.
Then
\[
\lambda_1(M)\leq\frac{(n-1+\theta)^2}{4}.
\]
Moreover, if $f$ be a positive function satisfying
\[
\Delta_\varphi f=-\lambda f
\]
for some constant $\lambda\geq 0$, then $f$ must satisfy the
gradient estimate
\[
|\nabla \ln f|^2\leq\frac{(n-1+\theta)^2}{2}-\lambda
+\sqrt{\frac{(n-1+\theta)^4}{4}-(n-1+\theta)^2\lambda}.
\]
\end{thm}

\vspace{0.5em}

After we have finished the proof of Theorem \ref{Rmain}, we saw that
Su-Zhang \cite{[SuZh]} and  Munteanu-Wang \cite{[MuWa2]} have
independently proved Theorem \ref{wth2}. We include it for completeness.
Our proof may be slight different from theirs. Here our approach
mainly follows the same spirit of proving Theorem 1.5 in \cite{[Wu2]}.

\section{Basic definitions and lemmas}\label{BaRe}
In this section, we will summarize some definitions and basic results on
a complete smooth metric measure space $(M^n,g,e^{-\varphi}dv)$.
For proofs and further details we refer to \cite{[Wang2]} and \cite{[Wu2]}.
\begin{definition}
Let $(M^n,g,e^{-\varphi}dv)$ be a complete smooth metric measure space.
A weighted Green's function $G_\varphi(x,y)$ is a function
defined on $(M\times M)\backslash\{(x,x)\}$ satisfying
$G_\varphi(x,y)=G_\varphi(y,x)$ and
$\Delta_{\varphi,y}G(x,y)=-\delta_{\varphi,x}(y)$ for all $x\neq y$,
where $\delta_{\varphi,x}(y)$ is defined by
\[
\int_M\psi(y)\delta_{\varphi,x}(y)d\mu=\psi(x)
\]
for every compactly supported function $\psi$.
\end{definition}
\begin{definition}
A complete smooth metric measure manifold $(M,g,e^{-\varphi}dv)$ is
said to be weighted non-parabolic if it admits a positive
weighted Green's function. Otherwise, it is said to be weighted parabolic.
\end{definition}
\begin{definition}
An end, $E$, with respect to a compact subset
 $\Omega\subset M$ is an unbounded connected component
of $M\setminus \Omega$. The number of ends with respect of
$\Omega$, denoted by $N_\Omega(M)$, is the number of unbounded connected
component of $M\setminus \Omega$.
\end{definition}
In general, when we say that $E$ is an end we mean that it is an end
with respect to some compact subset $\Omega$. In particular, its boundary
$\partial E$ is given by $\partial\Omega\cap \bar{E}$.

\begin{definition}
An end $E$ is said to be weighted non-parabolic if it admits a
positive weighted Green's function with Neumann boundary
condition on $\partial E$. Otherwise, it is said to be
weighted parabolic.
\end{definition}
Following similar arguments of Li-Tam \cite{[Li-Tam]}, we can
verify that a complete measure manifold is weighted non-parabolic
if and only if it has a weighted non-parabolic end. Of course,
it is possible for a weighted non-parabolic measure manifold to
have many weighted parabolic ends.

As similar discussions in Li-Tam \cite{[Li-Tam]}, we know
that an end $E$ with respect to the compact set $B_p(R_0)$ is
weighted non-parabolic if and only if there exists a sequence of
positive weighted harmonic functions $f_i$, defined on
$E_p(R_i)=E\cap B_p(R_i)$ for $R_0<R_1<R_2<\cdots \to\infty$,
satisfying $f_i=1$ on $\partial E$ and $f_i=0$ on
$\partial B_p(R_i)\cap E$, where $f_i$ is called
to be a barrier function of $E$. Moreover the sequence $f_i$
converges uniformly on compact subsets of $E\cup\partial E$ to
a minimal barrier function $f$, and $f$ has finite weighted
Dirichlet integral on $E$.

If $E$ be a weighted parabolic end with respect to $B_p(R_0)$,
then there exists a sequence of positive weighted harmonic
functions $g_i$, defined on $E_p(R_i)$ and the sequence of
constants $C_i\to\infty$, satisfying $g_i=0$
on $\partial E$ and $g_i=C_i$ on $\partial B_p(R_i)\cap E$.
Moreover the sequence $g_i$ converges uniformly on compact
subsets of $E\cap\partial E$ to a positive weighted function
$g$, satisfying $g=0$ on $\partial E$ and
$\sup_{y\in E}g=\infty$.

Now we give a decay estimate for weighted harmonic functions
on a weighted non-parabolic end of a smooth metric measure space,
which is a slight generalization of Lemma 1.1 in \cite{[Li-Wang1]}.
\begin{lemma}\label{lem0}
Let $(M^n,g,e^{-\varphi}dv)$ be a complete smooth metric measure space
of dimension $n$ with $\lambda_1(M)>0$. Assume that $E$ is a
weighted non-parabolic end of $M$. Then any weighted
harmonic function $f$ on $E$ satisfies the decay estimate
\[
\int_{E(R+1)\setminus E(R)}f^2d\mu\leq C\exp(-2\sqrt{\lambda_1(M)}R)
\]
for some constant $C>0$ depending on $f$, $\lambda_1(M)$ and $n$,
where $B_p(R)$ denotes a geodesic ball centered at some fixed
point $p\in M$ with radius $R>0$, and $E(R)=B_p(R)\cap E$.
\end{lemma}

The following lemma is an characterization for an end by its weighted
volume. Here we let $E$ be an end of $M$, and let $V_\varphi(E)$ be the
simply weighted volume of end $E$. We denote the weighted volume of
the set $E(R)$ by $V_\varphi(E(R))$.
\begin{lemma}\label{lem1}
Let $(M^n,g,e^{-\varphi}dv)$ be a complete smooth metric measure space
of dimension $n\geq 2$, satisfying $|\nabla\varphi|\leq\theta$ for some
positive constant $\theta$.  Assume that
$\lambda_1(M)\geq\frac{(n-1+\theta)^2}{4}$.
\begin{enumerate}
  \item  If $E$ is a weighted parabolic end, then
\[
V_\varphi(E)-V_\varphi(E(R))\leq C\exp(-(n-1+\theta)R)
\]
for some constant $C>0$ depending on $E$, where $R>0$ is large enough.

\item If $E$ is a weighted non-parabolic end, then
\[
V_\varphi(E(R))\geq C\exp((n-1+\theta)R)
\]
for some constant $C>0$ depending on $E$, where $R>0$ is large enough.
\end{enumerate}
\end{lemma}
\begin{remark}
Lemma \ref{lem1} can be viewed as a refined version of
eigenvalue estimate in Theorem \ref{wthm}. Indeed, if
$Ric_{\infty}\geq-(n-1)$ and $|\nabla\varphi|\leq\theta$
for constant $\theta>0$, then the weighted Bishop volume
comparison theorem (Theorem 1.2 in \cite{WeiWy}) asserts that
\[
V_\varphi(B_p(R))\leq C\exp(\theta R)\cdot
V_{\mathbb{H}^n}(B_p(R))\leq C\exp((n-1+\theta)R).
\]
Combining this  and Lemma \ref{lem1}, we conclude that
\[
\lambda_1(M)\leq\frac{(n-1+\theta)^2}{4},
\]
as asserted in Theorem \ref{wthm}. Regarding this, we think that
this trick is an effective method for the first eigenvalue upper
estimate on complete manifolds.
\end{remark}

On the other hand, if $Ric_{\infty}\geq-(n-1)$ and
$|\nabla\varphi|\leq\theta$ for constant $\theta>0$,
then the weighted Bishop volume comparison theorem,
asserts that for any $x\in M$ and $R_1<R_2$,
\[
\frac{V_\varphi(B_x(R_2))}{V_\varphi(B_x(R_1))}
\leq e^{\theta R_2}\frac{V_{\mathbb{H}^n}(B(R_2))}{V_{\mathbb{H}^n}(B(R_1))},
\]
where
$V_\varphi(B_x(R))=\int_{B_x(R)}e^{-\varphi}dv(g)$ denotes the weighted
volume of the geodesic ball $B_x(R)$, and $V_{\mathbb{H}^n}(B(R))$
denotes the volume of a geodesic ball of radius $R$ in the $n$-dimensional
hyperbolic space form $\mathbb{H}^n$ with constant curvature $-1$.
In particular, if we let $x=p$, $R_1=0$ and $R_2=R$, then
\begin{equation}\label{bijiao1}
V_\varphi(B_p(R))\leq C\exp((n-1+\theta)R)
\end{equation}
for sufficiently large $R$. If we let $x\in\partial B_p(R)$, $R_1=1$
and $R_2=R+1$, then
\begin{equation}
\begin{aligned}\label{bijiao2}
V_\varphi(B_x(1))&\geq CV_\varphi(B_x(R+1))\exp(-\theta)\exp(-(n-1+\theta)R)\\
&\geq C\exp(-\theta)V_\varphi(B_p(1))\exp(-(n-1+\theta)R).
\end{aligned}
\end{equation}
Combining \eqref{bijiao1}, \eqref{bijiao2} and Lemma \ref{lem1} immediately yields
that
\begin{corollary}\label{lem2}
Let $(M^n,g,e^{-\varphi}dv)$ be a complete smooth metric measure space
of dimension $n\geq 2$, satisfying $|\nabla\varphi|\leq\theta$ for some
positive constant $\theta$.  Assume that
$\lambda_1(M)\geq\frac{(n-1+\theta)^2}{4}$.
\begin{enumerate}
  \item  If $E$ is a weighted-parabolic end, then E must have exponential
weighted volume decay given by
\[
C_4\exp(-(n-1+\theta)R)\leq V_\varphi(E)-V_\varphi(E(R))\leq C_1\exp(-(n-1+\theta)R)
\]
for some constant $C_1\geq C_4>0$ depending on $E$ and $\theta$,
where $R>0$ is large enough.

\item If $E$ is a weighted-non-parabolic end, then $E$ must
have exponential volume growth given by
\[
 C_3\exp((n-1+\theta)R)\geq V_\varphi(E)\geq C_2\exp((n-1+\theta)R)
\]
for some constant $C_3\geq C_2>0$ depending on $E$, where $R>0$ is large enough.
\end{enumerate}
\end{corollary}

Besides the above properties, we can also confirm that if the spectrum of
the weighted Laplacian has a positive lower bound, then an end is
weighted non-parabolic if and only if its weighted volume is infinite
(see also \cite{[Wang2]}).


\section{Generalized Bochner formula}\label{gBocher}
In this section, we will give an improved version of the Bochner formula
for weighted harmonic functions, which is a mild generalization of the
case of harmonic functions due to Yau \cite{[Yau]}.
\begin{theorem}\label{wthbc}
Let $(M^n,g,e^{-\varphi}dv)$ be a complete smooth metric measure space
of dimension $n\geq2$ with $Ric_\infty\geq-(n-1)$. Assume that
$|\nabla\varphi|\leq\theta$ for some positive constant $\theta$ and
$f$ is a weighted harmonic function. For any constants $k>n$ and $p$,
\[
\Delta_\varphi|\nabla f|^p\geq\frac 1p
\left(\frac {k}{k-1}+p-2\right)|\nabla f|^{-p}(\nabla |\nabla f|^p)^2
-p\left(\frac{\theta^2}{k-n}+n-1\right)|\nabla f|^p.
\]
In particular, if $p=\frac{k-2}{k-1}$, then
\begin{equation}\label{guan1}
\Delta_\varphi|\nabla f|^p\geq-\frac{k-2}{k-1}
\left(\frac{\theta^2}{k-n}+n-1\right)|\nabla f|^p.
\end{equation}
\end{theorem}

\begin{proof}[Proof of Theorem \ref{wthbc}]
The proof is similar to the proof of Lemma 7.2 in \cite{[PLi]}. Choose
a local orthogonal frame $\{e_1, e_2,...,e_n\}$ near any
such given point so that at the given point $\nabla f=|\nabla
f|e_1$. Since $f$ is a weighted harmonic function and $Ric_\infty\geq-(n-1)$,
we compute that
\begin{equation}
\begin{aligned}\label{jisuan2}
\Delta_\varphi|\nabla f|^2&=\Delta|\nabla f|^2-\langle\nabla\varphi,\nabla|\nabla f|^2\rangle\\
&=2f_{ij}^2+2(R_{ij}+\nabla^2\varphi)f_if_j\\
&\geq 2f^2_{ij}-2(n-1)|\nabla f|^2.
\end{aligned}
\end{equation}
Notice that
\begin{equation}
\begin{aligned}\label{jisuan3}
\left|\nabla\left|\nabla f\right|^2\right|^2=4\sum^n_{j=1}
\left(\sum^n_{i=1}f_if_{ij}\right)^2
=4f^2_1\cdot\sum^n_{i=1}f^2_{1i}=4\left|\nabla
f\right|^2\cdot\sum^n_{i=1}f^2_{1i}.
\end{aligned}
\end{equation}
and
\begin{equation}
\begin{aligned}\label{holderin}
f_{ij}^2&\geq f^2_{11}+2\sum^n_{\alpha=2}f^2_{1\alpha}
+\sum^n_{\alpha=2}f^2_{\alpha\alpha}\\
&\geq f^2_{11}+2\sum^n_{\alpha=2}f^2_{1\alpha}
+\frac{1}{n-1}\left(\sum^n_{\alpha=2}f_{\alpha\alpha}\right)^2\\
&=f^2_{11}+2\sum^n_{\alpha=2}f^2_{1\alpha}
+\frac{1}{n-1}\left(\Delta f-f_{11}\right)^2\\
&=f^2_{11}+2\sum^n_{\alpha=2}f^2_{1\alpha}
+\frac{1}{n-1}\left(\varphi_if_i-f_{11}\right)^2\\
&\geq f^2_{11}+2\sum^n_{\alpha=2}f^2_{1\alpha}
+\frac{1}{n-1}\left[\frac{f_{11}^2}{1+\frac{k-n}{n-1}}
-\frac{(\varphi_if_i)^2}{\frac{k-n}{n-1}}\right]\\
&\geq \frac{k}{k-1}\sum^n_{j=1}f^2_{1j}
-\frac{\theta^2}{k-n}|\nabla f|^2
\end{aligned}
\end{equation}
for any constant $k>n$, where we used the following
inequality:
\[
(a+b)^2\geq\frac{a^2}{1+\delta}-\frac{b^2}{\delta}
\]
for any $\delta>0$.
Therefore
\begin{equation*}
\begin{aligned}
\Delta_\varphi|\nabla f|^2&\geq
\frac{k}{2(k-1)}|\nabla f|^{-2}|\nabla|\nabla f|^2|^2
-2\left(\frac{\theta^2}{k-n}+n-1\right)|\nabla f|^2\\
&=\frac{2k}{k-1}|\nabla|\nabla f||^2
-2\left(\frac{\theta^2}{k-n}+n-1\right)|\nabla f|^2.
\end{aligned}
\end{equation*}
Since
\[
\nabla|\nabla f|^p=\frac p2|\nabla f|^{p-2}\nabla|\nabla f|^2\quad
\mathrm{and}\quad\nabla|\nabla f|^2=2\nabla |\nabla f|\cdot\nabla|\nabla f|,
\]
we have that
\begin{equation*}
\begin{aligned}
\Delta_\varphi|\nabla f|^p&=\frac 1p(p-2)|\nabla f|^{-p}(\nabla |\nabla f|^p)^2
+\frac p2|\nabla f|^{p-2}\Delta_\varphi|\nabla f|^2\\
&\geq\frac 1p(p-2)|\nabla f|^{-p}(\nabla |\nabla f|^p)^2\\
&+\frac p2|\nabla f|^{p-2}\left[\frac{2k|\nabla f|^{2-2p}}{p^2(k-1)}(\nabla|\nabla f|^p)^2
-2\left(\frac{\theta^2}{k-n}+n-1\right)|\nabla f|^2\right]\\
&=\frac 1p\left(\frac {k}{k-1}+p-2\right)|\nabla f|^{-p}(\nabla|\nabla f|^p)^2
-p\left(\frac{\theta^2}{k-n}+n-1\right)|\nabla f|^p.
\end{aligned}
\end{equation*}
Letting $p=\frac{k-2}{k-1}$ in the above inequality yields \eqref{guan1}.
\end{proof}

\section{Proof of Theorem \ref{Rmain}}
In this section, we will prove Theorem \ref{Rmain} stated in introduction.
Here we mainly follow the arguments of Li-Wang \cite{[Li-Wang1]}.
\begin{proof}[Proof of Theorem \ref{Rmain}]
Assume that $M$ has at least two weighted non-parabolic ends.
By the construction described in Section \ref{BaRe},
there exists a nonconstant weighted harmonic function $f$ with
finite weighted Dirichlet integral on $M$.
Let $g$ be the function defined by
\[
g:=|\nabla f|^{\frac{k-2}{k-1}},
\]
where $k>n$ is a constant. By Theorem \ref{wthbc}, we have
\begin{equation}\label{subharm}
\Delta_\varphi g\geq-\frac{k-2}{k-1}
\left(\frac{\theta^2}{k-n}+n-1\right)g.
\end{equation}
In the following we \emph{claim} that the function $g$ satisfies
\begin{equation}\label{clam}
\int_{B_p(2R)\setminus B_p(R)}g^2d\mu\to 0
\end{equation}
for $R$ large enough. To prove this claim, firstly, we see that
\begin{equation}
\begin{aligned}\label{jisu0}
\int_{B_p(2R)\setminus B_p(R)}g^2d\mu\leq&
\left(\int_{B_p(2R)\setminus B_p(R)}\exp\left(2\sqrt{\lambda_1}r\right)|\nabla f|^2d\mu\right)^{\frac{k-2}{k-1}}\\
&\times \left(\int_{B_p(2R)\setminus B_p(R)}\exp\left(-2(k-2)\sqrt{\lambda_1}r\right)
d\mu\right)^{\frac{1}{k-1}},
\end{aligned}
\end{equation}
where we used the Holder inequality. In the following, we shall
estimate the right hand side of \eqref{jisu0}.
On one hand, since the weighted volume comparison theorem
\cite{{WeiWy}} asserts that
\[
A_{\varphi}(B_p(r))\leq \exp(\theta r)\cdot A_{\mathbb{H}^n}(B_p(r)),
\]
using this, we have
\begin{equation}
\begin{aligned}\label{guji1}
\int_{B_p(2R)\setminus B_p(R)}&\exp\left(-2(k-2)\sqrt{\lambda_1}r\right)d\mu\\
&\leq C\int^{2R}_R\exp\left(-2(k-2)\sqrt{\lambda_1}r\right)\exp(\theta r)\cdot\exp((n-1)r)dr\\
&=C\int^{2R}_R\exp\left[\left(n-1+\theta-2(k-2)\sqrt{\lambda_1}\right)r\right]dr
\end{aligned}
\end{equation}
Noticing
\[
\lambda_1(M)\geq\frac{k-2}{k-1}
\left(\frac{\theta^2}{k-n}+n-1\right),
\]
and $k>n\geq 3$, by Appendix \ref{app2}, we conclude that
\[
n-1+\theta-2(k-2)\sqrt{\lambda_1}<0
\]
and therefore the right hand side of \eqref{guji1} exponentially
decays to $0$. On the other hand, since we have assume that
$M$ has at least two weighted non-parabolic ends,
by Lemma 2.4 in \cite{[Wang2]}, we have the following decay estimate
\[
\int_{B_p(2R)\setminus B_p(R)}\exp\left(2\sqrt{\lambda_1}r\right)
|\nabla f|^2d\mu\leq CR.
\]
Combining this and \eqref{guji1}, we confirm our claim \eqref{clam}.

Next step, we consider a smooth compactly supported
function $\psi$ on $M$, satisfying
\begin{equation}
\psi(x)=\left\{
\begin{aligned}\label{cutfun}
1 \quad\quad& x\in B_p(R)\\
|\nabla \psi|\leq \frac CR \quad\quad& x\in B_p(2R)\setminus B_p(R)\\
0 \quad\quad& x\not\in B_p(2R).
\end{aligned}
\right.
\end{equation}
Then we have
\begin{equation}\label{inden1}
\int_M|\nabla(\psi g)|^2d\mu=\int_M|\nabla\psi|^2 g^2d\mu+\int_M\psi^2|\nabla g|^2d\mu
+2\int_M\psi g \nabla\psi \nabla gd\mu.
\end{equation}
Since
\[
2\int_M\psi g \nabla\psi \nabla gd\mu
=-\int_M\psi^2 |\nabla g|^2d\mu-\int_M\psi^2g\, \Delta_\varphi gd\mu,
\]
by \eqref{inden1}, we have
\begin{equation}
\begin{aligned}\label{jisua1}
\int_M|\nabla(\psi g)|^2d\mu
&=\int_M|\nabla\psi|^2 g^2d\mu-\int_M\psi^2g\, \Delta_\varphi gd\mu\\
&=\int_M|\nabla\psi|^2 g^2d\mu+\frac{k-2}{k-1}
\left(\frac{\theta^2}{k-n}+n-1\right)\int_M\psi^2g^2\\
&\quad-\int_M\psi^2g\left[\frac{k-2}{k-1}
\left(\frac{\theta^2}{k-n}+n-1\right)g+\Delta_\varphi g\right]d\mu.
\end{aligned}
\end{equation}
Since
\[
\lambda_1(M)\geq\frac{k-2}{k-1}
\left(\frac{\theta^2}{k-n}+n-1\right),
\]
the variational property of $\lambda_1(M)$ gives
\[
\frac{k-2}{k-1}
\left(\frac{\theta^2}{k-n}+n-1\right)
\int_M\psi^2g^2d\mu\leq\int_M|\nabla(\psi g)|^2d\mu.
\]
Substituting this into  \eqref{jisua1} yields
\begin{equation}\label{jisua2}
\int_M\psi^2g\left[\frac{k-2}{k-1}
\left(\frac{\theta^2}{k-n}+n-1\right)g+\Delta_\varphi g\right]d\mu
\leq \int_M|\nabla\psi|^2 g^2d\mu.
\end{equation}
We also note that the cut-function $\psi$ is the form of \eqref{cutfun},
and the right hand side of \eqref{jisua2} can be estimated by
\[
\int_M|\nabla\psi|^2 g^2d\mu\leq \frac{C}{R^2}
\int_{B_p(2R)\setminus B_p(R)}g^2d\mu.
\]
By our claim \eqref{clam}, we know the right hand side of integral tends to
$0$ as $R\to\infty$. Combining this with \eqref{subharm} we conclude that
$g$ either be to $0$ or it satisfies
\[
\Delta_\varphi g=-\frac{k-2}{k-1}
\left(\frac{\theta^2}{k-n}+n-1\right)g.
\]
By the construction of the weighted harmonic function $f$,
if $M$ has more than one weighted non-parabolic end, then function $f$ must be
non-constant, hence $g\neq 0$. So all the inequalities used to derive
\eqref{subharm} are equalities. In particular, from \eqref{holderin}
we conclude that
\[
f_{ij}=0 ,\quad i\neq j\quad \quad  \mathrm{and} \quad \quad
 f_{\alpha\alpha}=\rho, \quad 2\leq\alpha\leq n.
\]
The equality case of \eqref{subharm} also implies that
\[
|\nabla \varphi|=\theta\quad\quad\mathrm{and}\quad\quad
f_{11}=\frac{k-1}{k-n}\langle\nabla f,\nabla \varphi\rangle.
\]
Since $\Delta_\varphi f=0$, we also have
\[
f_{11}+(n-1)\rho=\langle\nabla f,\nabla \varphi\rangle.
\]
From above equalities, we derive that
\begin{equation}\label{guankx0}
f_{11}=-(k-1)\rho\quad\quad\mathrm{and}\quad\quad
\langle\nabla f,\nabla \varphi\rangle=-(k-n)\rho.
\end{equation}
Using the fact that $f_{1\alpha}=0$ for all $\alpha\neq 1$,
we conclude that $|\nabla f|$ is identically constant along the level set of
$f$. In particular, the level sets of $|\nabla f|$ and $f$ coincide. Moreover,
\[
\rho\delta_{\alpha\beta}=f_{\alpha\beta}=\mathrm{II}_{\alpha\beta} f_1,
\]
where $(\mathrm{II}_{\alpha\beta})$ is the second fundamental form of the level
sets of $f$. From this, we also have
\begin{equation}\label{guankx1}
f_{11}=-\frac{k-1}{n-1}Hf_1,
\end{equation}
where $H$ denotes the mean curvature of the level sets of $f$.
Applying the same computation to the function $g$, we have
\begin{equation}\label{guankx2}
-\frac{k-2}{k-1}\left(\frac{\theta^2}{k-n}+n-1\right)g
=\Delta_\varphi g=g_{11}+Hg_1-\langle\nabla g,\nabla \varphi\rangle.
\end{equation}
Here since $g=|\nabla f|^{\frac{k-2}{k-1}}$, we have
\[
g_1=\left(|\nabla f|^{\frac{k-2}{k-1}}\right)_1
=\frac{k-2}{k-1}|\nabla f|^{-\frac{k}{k-1}}f_jf_{j1}
=\frac{k-2}{k-1}f_1^{-\frac{1}{k-1}}f_{11}.
\]
Combining this with \eqref{guankx1} yields
\begin{equation}\label{dengshg}
H=-\frac{n-1}{k-2}g_1g^{-1}.
\end{equation}
Meanwhile, by the above equality, \eqref{guankx0} and \eqref{guankx1}
we have
\[
\varphi_1=\frac{k-n}{k-2}g_1g^{-1}.
\]
Plugging the above two relations into \eqref{guankx2}, we conclude that
\[
g_{11}-\frac{k-1}{k-2}g^2_1g^{-1}
+\frac{k-2}{k-1}\left(\frac{\theta^2}{k-n}+n-1\right)g=0.
\]
Letting
\[
u:=g^{-\frac{1}{k-2}}=|\nabla f|^{-\frac{1}{k-1}},
\]
the above equation becomes
\[
u_{11}-\frac{1}{k-1}\left(\frac{\theta^2}{k-n}+n-1\right)u=0.
\]
This equation can be regarded as as an ordinary differential equation
along the integral curve generated by the vector field $e_1$. Hence
\[
u(t)=A\exp\left[\sqrt{K(k,n,\theta)}\,\,t\right]
+B\exp\left[-\sqrt{K(k,n,\theta)}\,\,t\right],
\]
where $A$ and $B$ are nonnegative constants, and $K(k,n,\theta):=\frac{1}{k-1}\left(\frac{\theta^2}{k-n}+n-1\right)$.

Since $M$ has at least two weighted non-parabolic ends by the assumption
of theorem, we claim that any fixed level set $N$ of $|\nabla f|$
must be compact. Indeed, The facts that $f$ has no critical points
and that the level set of $f$ coincides with the level set of
$|\nabla f|$ lead to the splitting of $M$ as a warped product
$\mathbb{R}\times N$. Here the manifold $N$ is necessarily
compact due to the fact that $M$ is assumed to have at least
two ends. Since
\[
\int_M|\nabla f|^2d\mu<+\infty,
\]
we conclude that $|\nabla f|$ must have an interior maximum, saying
$\nabla f=1$. Now we fix $N=\{|\nabla f|=1\}$ and then $u$ must have its
minimum along $N$. By reparameterizing, we assume $N$ is given by $t=0$.
So
\[
0=u'(0)=A-B  \quad\mathrm{and}\quad 1=u(0)=A+B.
\]
This gives us
\[
u(t)=\cosh\left[\sqrt{K(k,n,\theta)}\,\,t\right]
\]
and
\[
g(t)=\cosh^{-(k-2)}\left[\sqrt{K(k,n,\theta)}\,\,t\right].
\]
Substituting this into \eqref{dengshg}, we have
\[
H(t)=(n-1)\sqrt{K(k,n,\theta)}
\tanh\left[\sqrt{K(k,n,\theta)}\,\,t\right]
\]
and
\[
h_{\alpha\beta}(t)=\sqrt{K(k,n,\theta)}
\tanh\left[\sqrt{K(k,n,\theta)}\,\,t\right]
\delta_{\alpha\beta}
\]
This implies that the metric on $M=\mathbb{R}\times N$ must be of the form
\eqref{biaoxing}.
\end{proof}

\section{Proof of Theorem \ref{wth2}}
We now follow the lines of \cite{[Wu2]} and prove Theorem \ref{wth2} given
in introduction. The proof method belongs to Li-Wang \cite{[Li-Wang2]}.
\begin{proof}[Proof of Theorem \ref{wth2}]
Since manifold $M$ satisfies the hypothesis of Theorem \ref{wth2},
Corollary \ref{Rcor} asserted that $M$ must have only one infinite weighted
volume end because the warped product with the metric given by
\[
ds_M^2=dt^2+\cosh tds_N^2
\]
has $\lambda_1(M)=n-2+\theta$. Now we assume that manifold $M$ has a
finite weighted volume end. Since $\lambda_1(M)>0$, $M$ must also
have an infinite weighted volume end. By choosing the compact set
$D$ appropriately, we may assume that $M\backslash D$ has one
infinite weighted volume, weighted non-parabolic end $E_1$ and
one finite weighted volume, weighted parabolic end $E_2$.

Following similar Li-Tam's arguments \cite{[Li-Tam]}, there exists
a positive weighted harmonic function $f$, satisfying
$\inf_{\partial E_1(R)}\to 0$  as $R\to\infty$,
$\sup_{\partial E_2(R)}\to \infty$  as $R\to\infty$  and
$f$ has finite weighted Dirichlet integral on $E_1$.
Meanwhile $f$ satisfies the following gradient estimate of
Theorem \ref{wthm}
\[
|\nabla f|^2\leq (n-1+\theta)^2f^2.
\]
Since function $f$ is weighted harmonic, we have
\begin{equation}\label{jisuan0}
\Delta_\varphi f^{1/2}=-\frac 14f^{-3/2}|\nabla f|^2
\leq-\frac{(n-1+\theta)^2}{4}f^{1/2}.
\end{equation}
If we let $h=f^{1/2}$, then for any nonnegative cut-off function $\psi$, we have
\begin{equation}\label{indent1}
\int_M|\nabla(\psi h)|^2d\mu=\int_M|\nabla\psi|^2 h^2d\mu+\int_M\psi^2|\nabla h|^2d\mu
+2\int_M\psi h \nabla\psi \nabla hd\mu.
\end{equation}
Since
\[
\int_M\psi h \nabla\psi \nabla hd\mu
=-\int_M\psi \nabla\psi h\nabla hd\mu
-\int_M\psi^2 |\nabla h|^2d\mu-\int_M\psi^2h\, \Delta_\varphi hd\mu,
\]
the integral equality \eqref{indent1} reduces to
\begin{equation}
\begin{aligned}\label{jisuan1}
\int_M|\nabla(\psi h)|^2d\mu
&=\int_M|\nabla\psi|^2 h^2d\mu-\int_M\psi^2h\, \Delta_\varphi hd\mu\\
&=\int_M|\nabla\psi|^2 h^2d\mu+\frac{(n-1+\theta)^2}{4}\int_M\psi^2h^2\\
&\quad-\int_M\psi^2h\left[\frac{(n-1+\theta)^2}{4}h+\Delta_\varphi h\right]d\mu.
\end{aligned}
\end{equation}
Since $\lambda_1(M)\geq\frac{(n-1+\theta)^2}{4}$, from the definition of
$\lambda_1(M)$, we have
\[
\frac{(n-1+\theta)^2}{4}\int_M\psi^2h^2d\mu\leq\int_M|\nabla(\psi h)|^2d\mu.
\]
Hence
\begin{equation}\label{jisuan2}
\int_M\psi^2h\left[\frac{(n-1+\theta)^2}{4}h+\Delta_\varphi h\right]d\mu
\leq \int_M|\nabla\psi|^2 h^2d\mu.
\end{equation}

Integrating the gradient estimate of Theorem \ref{wthm} along geodesics,
we know that $f$ must satisfy the growth estimate
\[
f(x)\leq C\exp((n-1+\theta)r(x)),
\]
where $r(x)$ is the geodesic distance from $x$ to a fixed point $p\in M$.
In particular, when restricted on the parabolic end $E_2$, together with
the volume estimate of Lemma \ref{lem1}, we conclude that
\begin{equation}\label{jisuan3}
\int_{E_2(R)}fd\mu
\leq CR.
\end{equation}
On the other hand, by Lemma \ref{lem0}, function $f$ must satisfy
the decay estimate on $E_1$
\[
\int_{E_1(R+1)\setminus E_1(R)}f^2d\mu\leq C\exp(-(n-1+\theta)R)
\]
for $R$ sufficiently large. By the Schwarz inequality, we have
\[
\int_{E_1(R+1)\setminus E_1(R)}fd\mu\leq C\exp\left(-\frac{n-1+\theta}{2}R\right)
V^{1/2}_{\varphi E_1}(R+1)
\]
where $V_{\varphi E_1}(r)$ denotes the weighted volume of $E_1(r)$.
Combining this with the volume estimate of Corollary \ref{lem2}, we
have that
\[
\int_{E_1(R+1)\setminus E_1(R)}fd\mu\leq C
\]
for some constant $C$ independent of $R$. In particular, we have
\[
\int_{E_1(R)}fd\mu\leq CR.
\]
Combining this with \eqref{jisuan3}, we conclude that
\begin{equation}\label{jisuan4}
\int_{B_p(R)}fd\mu \leq CR.
\end{equation}
Now we define the cut-off function $\psi$
on $M$ in \eqref{jisuan2} by
\begin{equation*}
\psi(x)=\left\{
\begin{aligned}
1 \quad\quad& x\in B_p(R)\\
\frac{2R-r}{R} \quad\quad& x\in B_p(2R)\setminus B_p(R)\\
0 \quad\quad& x\not\in B_p(2R).
\end{aligned}
\right.
\end{equation*}
Hence the right hand side of \eqref{jisuan2} is given by
\[
\int_M|\nabla\psi|^2 h^2d\mu=R^{-2}\int_{B_p(2R)\setminus B_p(R)} h^2d\mu
\]
and \eqref{jisuan4} implies
\[
\int_M|\nabla \psi|^2h^2d\mu \to 0
\]
as $R\to\infty$. Therefore we obtain
\[
\Delta_\varphi h=-\frac{(n-1+\theta)^2}{4}h
\]
and inequality \eqref{jisuan0} used in the above argument is an equality. In
particular,
\[
|\nabla f|=(n-1+\theta)f
\]
and hence
\begin{equation}\label{equ}
|\nabla (\ln f)|^2=(n-1+\theta)^2.
\end{equation}
Therefore inequalities used to prove the gradient estimate of
Theorem \ref{wthm} are all equalities. More precisely, we must have
equality (2.8) in \cite{[Wu1]} since
\[
\Delta_\varphi|\nabla(\ln f)|^2=\Delta|\nabla(\ln f)|^2
-\nabla \varphi\cdot\nabla|\nabla(\ln f)|^2=\nabla|\nabla(\ln f)|^2=0.
\]
Furthermore, inequalities used to derive (2.8) in \cite{[Wu1]} must all be
equalities. More specifically, equality (2.5) in \cite{[Wu1]} implies
\[
(\ln f)_{1j}=0
\]
for all $1\leq j\leq n$; whereas the equality of estimate term
$h^2_{ij}$ in \cite{[Wu1]} gives
\begin{equation}\label{equality2}
\langle\nabla \varphi,\nabla \ln f\rangle=(n-1+\theta)\theta
\end{equation}
and
\[
(\ln f)_{\alpha\beta}=-(n-1+\theta)\delta_{\alpha\beta}
\]
for all $2\leq \alpha,\beta\leq n$. Since $e_1$ is the unit normal to
the level set of $\ln f$, the second fundamental form $\mathrm{II}$ of
the level set is given by
\begin{equation*}
\begin{aligned}
\mathrm{II}_{\alpha\beta}&=\frac{(\ln f)_{\alpha\beta}}{(\ln f)_1}\\
&=\frac{-(n-1+\theta)\delta_{\alpha\beta}}{n-1+\theta}\\
&=-\delta_{\alpha\beta}.
\end{aligned}
\end{equation*}
Moreover, \eqref{equ} implies that if we set $t=\frac{\ln f}{n-1+\theta}$,
then $t$ must be the distance function between the level sets of $f$,
hence also for $\ln f$. Since
$\mathrm{II}_{\alpha\beta}=(-\delta_{\alpha\beta})$,
this implies that the metric on $M$ can be written as
\[
ds_M^2=dt^2+\exp(-2t)ds_N^2.
\]
By \eqref{equality2}, we also have
\[
\varphi(t,x)=\varphi(0,x)+\theta t,
\]
where $(t,x)\in\mathbb{R}\times N$. Since we assume that the
manifold $M$ has two ends, $N$ must be compact. A direct computation
shows that the condition $Ric_M\geq-(n-1)$
implies that  $Ric_N\geq 0$. This completes the proof of theorem.
\end{proof}

\section{Appendix}\label{app2}
In this part we will prove the following fact:
\begin{proposition}
if $k>n\geq 3$ and
\begin{equation}\label{eqass}
\lambda_1(M)\geq\frac{k-2}{k-1}
\left(\frac{\theta^2}{k-n}+n-1\right),
\end{equation}
then we have
\[
n-1+\theta-2(k-2)\sqrt{\lambda_1}<0.
\]
\end{proposition}

\begin{proof}
Indeed, we only need to confirm that
\[
2(k-2)\sqrt{\lambda_1}>n-1+\theta.
\]
That is,
\[
4(k-2)^2\lambda_1>(n-1+\theta)^2.
\]
We also notice that \eqref{eqass}. Hence if we can prove
\[
4(k-2)^2\frac{k-2}{k-1}
\left(\frac{\theta^2}{k-n}+n-1\right)>(n-1+\theta)^2,
\]
then the desired conclusion follows. For the above inequality, rearranging
terms gives
\[
\left[\frac{4(k-2)^3}{(k-1)(k-n)}-1\right]\theta^2-2(n-1)\theta
+\left[\frac{4(k-2)^3(n-1)}{k-1}-(n-1)^2\right]>0.
\]
This is a quadratic inequality in $\theta$. We assert that this
inequality is always true. Because when $k>n\geq 3$, we have
\[
\frac{4(k-2)^3}{(k-1)(k-n)}-1>0
\]
and
\begin{equation*}
\begin{aligned}
4(n-1)^2&-4\left[\frac{4(k-2)^3}{(k-1)(k-n)}-1\right]
\left[\frac{4(k-2)^3(n-1)}{k-1}-(n-1)^2\right]\\
&=\frac{16(k-2)^3(n-1)}{(k-1)^2(k-n)}\left[-4(k-2)^3+(k-1)^2\right]\\
&<0.
\end{aligned}
\end{equation*}
This proves the proposition.
\end{proof}

\bibliographystyle{amsplain}

\begin{thebibliography}{30}
\bibitem{[BE]}D. Bakry, M. \'{E}mery, Diffusion hypercontractivitives,
in: S\'{e}minaire de Probabilit\'{e}s XIX, 1983/1984, in: Lecture
Notes in Math., vol. 1123, Springer-Verlag, Berlin, 1985, pp.
177-206.

\bibitem{[BQ1]}D. Bakry, Z.-M. Qian, Some new results on eigenvectors via
dimension, diameter and Ricci curvature, Adv. Math. 155 (2000)
98-153.

\bibitem{[BQ2]}D. Bakry, Z.-M. Qian, Volume comparison theorems without Jacobi
fields, in: Current Trends in Potential Theory, in: Theta Ser. Adv.
Math., vol. 4, Theta, Bucharest, 2005, pp. 115-122.

\bibitem{[Cao1]}H.-D. Cao, Recent progress on Ricci solitons. In Recent Advances in
Geometric Analysis, volume 11 of Advanced Lectures in Mathematics
(ALM). International Press, 2009.

\bibitem{[Che-Gro1]} J. Cheeger, D. Gromoll, On the structure of complete
manifolds of nonnegative curvature, Ann. Math., 92 (1972) 413-443.

\bibitem{[Che-Gro2]} J. Cheeger, D. Gromoll, The splitting theorem for
manifolds of nonnegative Ricci curvature, J. Diff. Geom., 6 (1971) 119-128.

\bibitem{[Cheng]}S.-Y. Cheng, Eigenvalue comparison theorems and its geometric
applications, Math. Z. 143 (1975) 289-297.

\bibitem{[FLZ]} F.-Q. Fang, X.-D. Li, Z.-L. Zhang, Two generalizations of
Cheeger-Gromoll splitting theorem via Bakry-\'{E}mery Ricci curvature,
Annales de l'Institut Fourier 59 (2009) 563-573.

\bibitem{[Hamilton]}R. S. Hamilton, The formation of singularities in the
Ricci flow, Surveys in Differential Geometry, 2: 7-136, 1995.

\bibitem{[PLi]}P. Li, Harmonic functions and applications to complete
manifolds, ``http://math.uci.edu/~pli/", preprint, 2004.

\bibitem{[Li-Tam]} P. Li, L.-F. Tam, Symmetric Green's functions on complete
manifolds, Amer. J. Math. 109 (1987) 1129-1154.

\bibitem{[Li-Wang1]}P. Li, J.-P. Wang, Complete manifolds with positive
spectrum, J. Diff. Geom. 58 (2001) 501-534.

\bibitem{[Li-Wang2]}P. Li, J.-P. Wang, Complete manifolds with positive
spectrum, II, J. Diff. Geom. 62 (2002) 143-162.

\bibitem{[LD]}X.-D. Li, Liouville theorems for symmetric diffusion operators on
complete Riemannian manifolds, J. Math. Pure. Appl. 84 (2005)
1295-1361.

\bibitem{[Lott]} J. Lott, Some geometric properties of the Bakry-\'{E}mery Ricci
tensor, Comment. Math. Helv. 78 (2003) 865-883.

\bibitem {[MuWa]} O. Munteanu, J. Wang, Smooth metric measure spaces with
nonnegative curvature, Comm. Anal. Geom. 19 (2011) 451-486.

\bibitem {[MuWa2]} O. Munteanu, J. Wang, Analysis of weighted Laplacian
and applications to Ricci solitons, to appear in Comm. Anal. Geom.

\bibitem {[SuZh]} Y. Su and H. Zhang, Rigidity of manifolds with Bakry-\'{E}mery
Ricci curvature bounded below, Geome. Dedicata, DOI 10.1007/s10711-011-9685-x.

\bibitem{[Wang1]} L.-F. Wang, The upper bound of the $L_{\mu}^2$ spectrum,
 Ann. Glob. Anal. Geom. 37(4) (2010) 393-402.

\bibitem{[Wang2]} L.-F. Wang,  A splitting theorem for the weighted measure,
 Ann. Glob. Anal. Geom. DOI 10.1007/s10455-011-9302-0.

\bibitem{WeiWy} G.-F. Wei, W. Wylie, Comparison geometry for the
Bakry-\'{E}mery Ricci tensor, J. Diff. Geom., 83 (2009), 377-405.

\bibitem{[Wu]} J.-Y. Wu, Upper bounds on the first eigenvalue for a diffusion
operator via Bakry-\'{E}mery Ricci curvature, J. Math. Anal. Appl.
361 (2010) 10-18.

\bibitem{[Wu1]} J.-Y. Wu, Upper bounds on the first eigenvalue for a diffusion
operator via Bakry-\'{E}mery Ricci curvature II, arXiv: math.DG/1010.4175.

\bibitem{[Wu2]} J.-Y. Wu, A note on the splitting theorem for the weighted measure

\bibitem{[Yau]} S.-T. Yau, Harmonic functions on complete Riemannian manifolds, Comm.
Pure Appl. Math. 28 (1975) 201-228.

\end{thebibliography}

\end{document}